\documentclass[11pt]{article}
\usepackage[T1]{fontenc}
\usepackage[toc,page]{appendix} 
\usepackage{amsmath}
\usepackage{amsfonts}
\usepackage{amssymb}
\usepackage{geometry}
\usepackage{fancyhdr}
\author{Vincent ~NOLOT\\ \vspace{5pt} \\
   Institut de Mathématiques de Bourgogne,\\
   Université de Bourgogne, 21078 Dijon, France.\\
   \texttt{vincent.nolot@u-bourgogne.fr}}
\title{Monge Problem on infinite dimensional Hilbert space endowed with suitable Gaussian measure}

\date{}

\def\R{\mathbb{R}}
\def\N{\mathbb{N}}

\def\P{\mathcal{P}}

\def\L{\mathcal{L}}

\def\V{\mathcal{V}}

\newcommand{\ee}{\varepsilon}

\def\fin{\hfill$\square$}

\newtheorem{theorem}{Theorem}[section]
\newtheorem{lemma}[theorem]{Lemma}       
\newtheorem{corollary}[theorem]{Corollary}
\newtheorem{proposition}[theorem]{Proposition}
\newtheorem{remark}[theorem]{Remark}

\newtheorem{definition}[theorem]{Definition}

\begin{document}

\maketitle

\abstract{In this paper we solve the Monge problem on infinite dimensional Hilbert space endowed with a suitable
Gaussian measure.}


\section{Introduction}
Our framework is an infinite dimensional Hilbert space $(H,|.|)$ endowed with its Borelian $\sigma-$algebra.
For $\rho_0$ and $\rho_1$ two Borel probability measures on $H$,
the \emph{Monge Problem} consists of finding a Borel map $T:H\longrightarrow H$ satisfying the constraint
$T_\#\rho_0(B):=\rho_0(T^{-1}(B)) = \rho_1(B)$ (for any Borel subset $B$ of $H$) and minimizing the quantity
$$\int_H c(x,T(x))d\rho_0(x),$$
where $c:H\times H \longrightarrow [0,\infty)$ is called \emph{cost function}.

\begin{theorem}\label{MongeHilbert}
Assume that $\rho_0$ and $\rho_1$ have finite relative entropy with respect to $\gamma$ where $\gamma$ satisfies
conditions of Theorem \ref{LebPoints}. Then the problem
\begin{equation}\label{MongeHilbertMap}
	\inf_{T_\#\rho_0=\rho_1} \int_H |x-T(x)|d\rho_0(x)
\end{equation}
has at least one solution $T:H \longrightarrow H$.
\end{theorem}

Monge Problem has been solved in infinite dimensional Hilbert spaces, when the cost
is $c(x,y)=|x-y|^p$ and $p>1$ (see e.g. \cite{AMB}). The case when $p$ is equal to $1$ is quite more tricky.
This is the object of our paper.

We are inspired from Champion and De Pascale in \cite{CHAMP}. The strategy for infinite dimensional case 
lies on the same powerful tool as in finite dimensional case: an essential
ingredient is the \emph{differentiation theorem for the measure of reference}. Unfortunately there is some
measure on Hilbert spaces for which this theorem is false. 
Nevertheless Tiser has proved in \cite{TISER} that for a suitable Gaussian measure on some Hilbert space, 
the differentiation theorem holds, namely:

\begin{theorem}\label{LebPoints}
Let $H$ be a separable Hilbert space and let $\gamma$ be a Gaussian measure with the following representation
of its covariance operator :
$$R(x)=\sum_i c_i(x,e_i)e_i,$$
where $(e_i)_{i}$ is an orthonormal system of $H$. Suppose that for $\alpha >5/2$ given we have
$c_{i+1}\leq c_i/i^\alpha$ for all $i$. Then
$$\lim_{r\rightarrow 0}\frac{1}{\gamma(B(x,r))}\int_{B(x,r)}|f-f(x)|d\gamma =0~~~~\mathrm{for~}\gamma-\mathrm{a.a.~x}\in H$$
for any $f\in L^p(H,\gamma)$ and $p>1$.
\end{theorem}

The set of $x\in H$ such that Theorem \ref{LebPoints} holds, is called the set of \emph{Lebesgue points} of $f$
and will be denoted by $Leb(f)$. Thus $\gamma(Leb(f))=1$. In the case of $f=1\!\!1_{A}$, we will call $x$ a \emph{Lebesgue
point of} $A$.

\vskip 2mm

\begin{remark}
In fact the Theorem \ref{LebPoints} is required only to get the Proposition \ref{HSupProp}. All other results
in this section are available without Lebesgue points.
\end{remark}

\vskip 2mm

\quad From now, $\gamma$ is the Gaussian measure defined on $H$ satisfying conditions of the previous Theorem \ref{LebPoints}. 
So that the \emph{differentiation theorem} holds over $(H,\gamma)$.
\vspace{10pt}

The classical way to find a solution of (\ref{MongeHilbertMap}) is to introduce the following Monge-Kantorovich problem :
\begin{equation}\label{MK1Hilbert}
	\min_{\Pi\in C(\rho_0,\rho_1)}\int_{H \times H}|x-y|d\Pi(x,y),
\end{equation}
where $C(\rho_0,\rho_1)$ is the set of \emph{coupling between} $\rho_0$ and $\rho_1$. 
The nonempty set of solutions (optimal couplings) of (\ref{MK1Hilbert})
will be denoted by $\mathcal{O}_1(\rho_0,\rho_1)$. Among these coupling, we shall show there is 
at least one which is carried by a graph of some map $T$ and therefore this map will be a solution of 
(\ref{MongeHilbertMap}).

Because the cost induced by the euclidian norm is not strictly convex, 
the set $\mathcal{O}_1(\rho_0,\rho_1)$ does not contain enough information to construct some map $T$.
Thus we need introduce an other problem, called \emph{second variational problem}, with a new cost to minimize over
the set of optimal couplings of (\ref{MK1Hilbert}):
\begin{equation}\label{MK2Hilbert}
	\min_{\Pi\in \mathcal{O}_1(\rho_0,\rho_1)}\int_{H \times H}\alpha(x-y)d\Pi(x,y),
\end{equation}
with 
$$\alpha(x-y):=\sqrt{1+|x-y|^2}.$$
This cost $\alpha$ is strictly convex and smooth. It turns out that it shall bring more information, 
namely in some sense the directions
that should take the optimal plan in order to be concentrated on a graph of some map.

We denote by $\mathcal{O}_2(\rho_0,\rho_1)$ the subset of $\mathcal{O}_1(\rho_0,\rho_1)$ containing
optimal couplings which minimize (\ref{MK2Hilbert}).
It is easy to see that $\alpha(x-y)\leq 1+|x-y|$ so that if (\ref{MK1Hilbert}) is finite for some coupling
then (\ref{MK2Hilbert}) is also finite, and
the set $\mathcal{O}_2(\rho_0,\rho_1)$ is a nonempty (by weak compacity) and a convex subset of
$C(\rho_0,\rho_1)$.

\vskip 2mm

For our purpose, we need to consider finite dimensional
approximations. Since the eigenvalues of the covariance matrice of $\gamma$ are $(c_i)_i$ 
we can identify $(H,\gamma)$ with $(l^2(c),\mu)$ where $\mu$ is the product of standard Gaussian measure on $\R$ and
$$l^2(c):=\{x\in \R^{\N},~\sum c_ix_i^2<\infty\}.$$
This latter space is separable, therefore we consider a sequence of maps $(\pi_n)_n$ such that
each $\pi_n$ projects $l^2(c)$ onto a $n-$dimensional euclidian space, and $\lim_n \pi_n = Id$ the identity map on $l^2(c)$.
In the sequel we make the abuse of notation $(H,\gamma)=(l^2(c),\mu)$.

\vskip 2mm

We denote by $\mathcal{A}_n$ the $\sigma-$algebra generated by $\pi_n$, and if (for $i=0,1$) $f_i$ is the density 
of $\rho_i$ w.r.t. $\mu$, we put $\hat{\rho}_i^n:=\mathbb{E}[f_i|\mathcal{A}_n]\gamma$.
For $x\in H$ and $n\in \mathbb{N}$, we denote by $x_n:=\pi_n(x)$.

We say that a coupling $\Pi\in C(\rho_0,\rho_1)$ satisfies the \emph{convexity property} if the relative entropy
is $1-$convex along geodesics $\rho_t:=((1-t)P_1+tP_2)_\#\Pi$, namely 
$$Ent_\gamma(\rho_t)\leq (1-t)Ent_\gamma(\rho_0)+tEnt_\gamma(\rho_1)-W^2(\rho_0,\rho_1),$$
holds for any $t\in (0,1)$.

Finally we are interested in the following set:
\begin{eqnarray*}
	 \overline{\mathcal{O}_2}(\rho_0,\rho_1):=\big\{\Pi \in \mathcal{O}_2(\rho_0,\rho_1),
\Pi~\mathrm{enjoys~the~\emph{convexity~property}}\big\}.
\end{eqnarray*}

The fact that $\overline{\mathcal{O}_2}(\rho_0,\rho_1)$ is non empty is the purpose of Theorem \ref{NonEmptyO2}. It
will play a key role in our approach because any coupling of $\overline{\mathcal{O}_2}(\rho_0,\rho_1)$
will bring us enough information to show that it is concentrated on a graph of some measurable map.

\vskip 2mm

Let us present how this paper is organized. In section 2 we establish the convexity of relative entropy 
(w.r.t. $\mu$) in
$(\mathcal{P}_1(H),W_1)$. In particular we obtain that if $\Pi \in \overline{\mathcal{O}_2}(\rho_0,\rho_1)$
then $\rho_t:=((1-t)P_1+tP_2)_\#\Pi$ (here $P_i$ designs the projection onto the $i-th$ component)
is absolutely continuous w.r.t. $\mu$ for any $t\in (0,1)$. This point will be
necessary through the next sections.
In section 3 we present different features of the support of element belonging to $\overline{\mathcal{O}_2}(\rho_0,\rho_1)$.
Proposition \ref{HSupProp}, which relies on Lemma \ref{SuppHilbLem1} is paramount for the method
used in the proof of Theorem \ref{MongeHilbert}.
Finally the last section is devoted to prove Theorem \ref{MongeHilbert}. The end contains many comments about the proof
and about open problems.

\section{Convexity of relative entropy in $(\mathcal{P}_1(H),W_1)$}

The following Proposition states that the relative entropy with respect to the Lebesgue measure on $\R^n$
is \emph{convex} along geodesics in $(\P_p(\R^n),W_p)$ whatever $p>1$. 
It is fundamental to get all other results of \emph{convexity} of relative entropy (when the reference measure
is absolutely continuous with respect to the Lebesgue measure).

\begin{proposition}\label{EntRnLeb}
	Let $c$ be a strictly convex and differentiable norm on $\mathbb{R}^n\backslash \{0\}$.
	If $p>1$ then for any $\rho_0,\rho_1\in D(Ent_{\L})$ and $\Pi$ optimal (for $c$) coupling between $\rho_0$
	and $\rho_1$, $\rho_t:=(T_t)_\#\Pi$ satisfies
	$$Ent_{\L}(\rho_t)\leq (1-t)Ent_{\L}(\rho_0)+tEnt_{\L}(\rho_1),~~~~\forall t\in [0,1].$$
\end{proposition}

{\bf Proof.}
	See for example \cite{VILL} (Chapter 5.).
\fin
\vskip 2mm

In order to extend our result in infinite dimensional spaces, we work with Gaussian measures as reference measures.
Let $\gamma_n$ be the standard Gaussian measure on $\R^n$.
We consider $\rho_0$ and $\rho_1$ two probability measures on $\R^n$ belonging to $D(Ent_\gamma)$.
For the Euclidian norm $|.|$ on $\R^n$, we introduce quantity \emph{inspired from} the so called Wasserstein distance:
$$\mathcal{W}_{\ee}(\rho_0,\rho_1):=\inf_{\Pi \in C(\rho_0,\rho_1)} \int_{\R^n\times \R^n}|x-y|+\ee \alpha(x-y)d\Pi(x,y),$$
where 
$$\alpha(x-y):=\left(1+|x-y|\right)^{1/2}.$$
Here $\alpha$ is strictly convex and differentiable function on $\R^n$.
We have the relation:
$$c_{\ee}(x-y):=|x-y|+\ee \alpha(x-y) \leq \ee+(1+\ee)|x-y|.$$
Beside $c_{\ee}$ is not a distance, neither is $\mathcal{W}_{\ee}$. 

\vskip 2mm

We recall that the $1-$Wasserstein distance in this situation is defined as
$$W_{1}(\rho_0,\rho_1):=\inf_{\Pi \in C(\rho_0,\rho_1)}\int_{\R^n\times \R^n}|x-y|\Pi(x,y).$$
Because of 
$\liminf_{\ee \rightarrow 0} \mathcal{W}_{\ee}(\rho_0,\rho_1) \geq W_{1}(\rho_0,\rho_1)$ 
we always fix $\ee$ small enough in such a way that,
$$\mathcal{W}_{\ee}(\rho_0,\rho_1)-\ee \geq W_{1}(\rho_0,\rho_1)-\ee>0.$$

\begin{proposition}\label{EntRnNormCee}
	If $\Pi$ is optimal for the cost $c_{\ee}$
	then for any $t\in (0,1)$ and:
	\begin{equation}\label{ConvEnt1FiniCee0}
	Ent_{\gamma_n}(\rho_t)\leq (1-t)Ent_{\gamma_n}(\rho_0)+tEnt_{\gamma_n}(\rho_1)-\frac{t(1-t)}{2(1+\ee)^2}
\left(\mathcal{W}_{\ee}(\rho_0,\rho_1)-\ee\right)^2.
	\end{equation}
	In particular if $\rho_0,\rho_1\in D(Ent_{\gamma_n})$ then also $\rho_t\in D(Ent_{\gamma_n})$ for any $t\in (0,1)$.
\end{proposition}

{\bf Proof.}
We can assume that $\rho_0,\rho_1\in D(Ent_{\gamma_n})$, otherwise the inequality is obvious.
Therefore since $\rho_0$ and $\rho_1$ be two probability measures absolutely continuous with respect to $\gamma_n$,
they are also absolutely continuous with respect to the Lebesgue measure $\mathcal{L}$. 
For $i=0,1$ let $d\rho_0=f_0d\mathcal{L}$ and $d\rho_1=f_1d\mathcal{L}$, then
the density of probability of $\rho_i$ with respect to $\gamma_n$ is $\frac{d\rho_i}{\gamma_n}=f_i(2\pi)^{\frac{d}{2}}
e^{\frac{|x|^2}{2}}$. Write:
\begin{eqnarray*}
	Ent_{\gamma_n}(\rho_i) &=&\int_{\mathbb{R}^n}f_i(x)(2\pi)^{\frac{d}{2}}e^{\frac{|x|^2}{2}}\log\left(f_i(x)
	(2\pi)^{\frac{d}{2}}e^{\frac{|x|_2^2}{2}}\right)d\gamma_n(x)\\
	&=&\int_{\R^n}	f_i(x)(2\pi)^{\frac{d}{2}}e^{\frac{|x|^2}{2}}\log(f_i(x))d\gamma_n(x)
		+\int_{\R^n}	f_i(x)(2\pi)^{\frac{d}{2}}e^{\frac{|x|^2}{2}}\log((2\pi)^{\frac{d}{2}})d\gamma_n(x)\\
	&~&
		+\int_{\R^n} f_i(x)(2\pi)^{\frac{d}{2}}e^{\frac{|x|^2}{2}}\frac{|x|^2}{2}d\gamma_n(x)\\
	&=& Ent_{\mathcal{L}}(\rho_i)+\V(\rho_i)+\frac{d}{2}\log(2\pi),
\end{eqnarray*}
where $\V(\rho_i):=\frac{1}{2}\int |x|_2^2d\rho_i(x)$.
	By $1-$convexity of the euclidian norm, it is	easy to see that:
	$$\mathcal{V}(\rho_t)\leq (1-t)\mathcal{V}(\rho_0)+t\mathcal{V}(\rho_1)-\frac{t(1-t)}{2}\int_{\R^n} |x-y|^2d\Pi(x,y).$$
	By Cauchy-Schwarz inequality, we get
	
\begin{equation}\label{ConvEntLem1}
\mathcal{V}(\rho_t)\leq (1-t)\mathcal{V}(\rho_0)+t\mathcal{V}(\rho_1)-\frac{t(1-t)}{2(1+\ee)^2}
\left(\mathcal{W}_{\ee,\|.\|}(\rho_0,\rho_1)-\ee\right)^2.
\end{equation}
then combining the Proposition \ref{EntRnLeb} with (\ref{ConvEntLem1}), we get the result taking the sum.
\fin

\vskip 2mm

Now we focus on our separable infinite dimensional Hilbert space $(H,\gamma)$.

\vskip 2mm

In order to apply these results above, we are interested in the following problem:
\begin{equation}\label{ProblemDerive}
\min_{\Pi\in C(\rho_0,\rho_1)} \int_{H\times H}|x-y|d\Pi(x,y)+\ee \int_{H\times H}\alpha(x-y)d\Pi(x,y), \tag{$P_\ee$}
\end{equation}
where $\alpha$ is defined as above,
$$\alpha(x-y):=\left(1+|x-y|\right)^{1/2}.$$
Here $|.|$ stands for the Hilbert norm on $H$.

\vskip 2mm

The following result extends the Proposition \ref{EntRnNormCee} to the infinite dimensional Hilbert space.

\begin{proposition}\label{HEnt0}
		Let $\Pi_\ee$ be a solution of $(P_\ee)$, being $w-$limit point
		of a sequence $(\Pi_n)_n$ with $\Pi_n \in C(\rho_0^n,\rho_1^n)$ optimal for $c_{\ee}$
		and satisfying (\ref{ConvEnt1FiniCee0}).
		If $\rho_t:=(T_t)_\#\Pi$ then for any $t\in (0,1)$, $\rho_t\in D(Ent_\gamma)$ and:
\begin{equation}\label{HConvEnt3Fini0}
Ent_\gamma(\rho_t)\leq (1-t)Ent_\gamma(\rho_0)+tEnt_\gamma(\rho_1)-\frac{t(1-t)}{2(1+\ee)^2}\left(\mathcal{W}_{\ee}(\rho_0,\rho_1)-\ee\right)^2.
\end{equation}
\end{proposition}

\vskip 2mm

{\bf Proof.}
	Fix $n\in \mathbb{N}$. Define $\rho_t^n:=(T_t)_\#\Pi_n$ for $t\in (0,1)$. 
	Because of all measures $\rho_i^n$ can be seen on probability measures over $H$,
$$	Ent_\gamma(\rho_t^n) = Ent_{\gamma_n}(\rho_t^n)~~~~\forall t\in [0,1],$$
	and we apply the Proposition \ref{EntRnNormCee} in the case of the Euclidian norm $|.|$, that is for all $t\in [0,1]$:
$$Ent_{\gamma}(\hat{\rho}_t^n)\leq (1-t)Ent_{\gamma}(\hat{\rho}_0^n)+tEnt_{\gamma}(\hat{\rho}_1^n)-\frac{t(1-t)}{2(1+\ee)^2}\left(\mathcal{W}_{\ee}(\rho_0^n,\rho_1^n)-\ee\right)^2.$$
	And $\mathcal{W}_{\ee,|.|}^2(\rho_0,\rho_1)\leq \liminf_n \mathcal{W}_{\ee,|.|}^2(\rho_0^n,\rho_1^n)$, therefore
	if $\delta >0$ is small enough so that $\mathcal{W}_{\ee,|.|}^2(\rho_0,\rho_1)-\ee-\delta>0$, 
	we can find $N\in \mathbb{N}$ such that:
	$$\mathcal{W}_{\ee,|.|}(\rho_0^n,\rho_1^n) +\delta \geq \mathcal{W}_{\ee,|.|}(\rho_0,\rho_1)~~~~\forall n \geq N.$$
	Jensen's inequality implies $Ent_\gamma(\rho_i^n)\leq Ent_\gamma(\rho_i)$ for $i=0,1$. Then for all $n\geq N$:
$$Ent_{\gamma}(\rho_t^n)\leq (1-t)Ent_{\gamma}(\rho_0)+tEnt_{\gamma}(\rho_1)-\frac{t(1-t)}{2(1+\ee)^2}\left(\mathcal{W}_{\ee}(\rho_0,\rho_1)-\ee-\delta\right)^2.$$
	Since $(\Pi_n)_n$ converges weakly to $\Pi$, it is the same for $(\rho_t^n)_n$ to $\rho_t$,
	and the compacity of the set $\{Ent_\gamma(.)\leq R\}$, and the lower semicontinuity of $Ent_\gamma(.)$ let us	to conclude:
$$Ent_{\gamma}(\rho_t)\leq (1-t)Ent_{\gamma}(\rho_0)+tEnt_{\gamma}(\rho_1)-\frac{t(1-t)}{2(1+\ee)^2}\left(\mathcal{W}_{\ee}(\rho_0,\rho_1)-\ee-\delta\right)^2.$$	
	Letting $\delta \rightarrow 0$, the result follows.	
\fin

\vskip 3mm

For the next Corollary, we deal with the true \emph{Wasserstein distance} $W_{1,|.|}$ on $\P(H)$. In this case for $\Pi \in \mathcal{O}_1(\rho_0,\rho_1)$ we can talk
about (constant speed) \emph{geodesics} for $\rho_t:=(T_t)_\#\Pi$, namely
$$W_1(\rho_t,\rho_s)=|t-s|W_1(\rho_0,\rho_1),~~~~\forall t\in [0,1].$$

\begin{corollary}\label{HEnt1}
	Let $\Pi \in C(\rho_0,\rho_1)$ be a $w-$limit point of $(\Pi_\ee)_\ee$ solutions of $(P_{\ee})$, and such
	that each $\Pi_\ee$ satisfies (\ref{HConvEnt3Fini0}). If $\rho_t:=(T_t)_\#\Pi$ then for any
	$t\in (0,1)$, $\rho_t\in D(Ent_\mu)$ and:
\begin{equation}\label{HConvEnt4Fini0}
Ent_\gamma(\rho_t)\leq (1-t)Ent_\gamma(\rho_0)+tEnt_\gamma(\rho_1)-\frac{t(1-t)}{2}W_{1}^2(\rho_0,\rho_1).
\end{equation}
\end{corollary}

In the literature, this proposition can be reformulated as: \emph{relative entropy is geodesically 
$1-$convex in $(\P(H),W_{1,|.|})$}.

\vskip 3mm

{\bf Proof.}
  Let $\rho_t^\ee:=((1-t)P_1+tP_2)_\#\Pi_\ee$. Thanks to the Proposition \ref{HEnt0}:
$$Ent_\gamma(\rho_t^\ee)\leq (1-t)Ent_\gamma(\rho_0)+tEnt_\gamma(\rho_1)-\frac{t(1-t)}{2(1+\ee)^2}\left(\mathcal{W}_{\ee,|.|}(\rho_0,\rho_1)-\ee\right)^2. $$
	Because $c_{\ee,|.|}$ converges to the Hilbert norm $|.|$ when
	$\ee$ goes to $0$, it turns out that $W_{1,|.|}(\rho_0,\rho_1)\leq \liminf_{\ee \rightarrow 0}
	\mathcal{W}_{\ee,|.|}(\rho_0,\rho_1)$.
	Arguing as in the proof above, for all $\delta>0$ and $\ee$ small enough:
$$Ent_\gamma(\rho_t)\leq (1-t)Ent_\mu(\rho_0)+tEnt_\gamma(\rho_1)-\frac{t(1-t)}{2(1+\ee)^2}\left(W_{1,|.|}(\rho_0,\rho_1)-\ee-\delta\right)^2. $$	
	Finally we let $\ee $ goes to $0$ and then $\delta$ goes to $0$.
\fin

\vskip 3mm

A particular case of application of the previous Proposition is the following :
if $\Pi\in \overline{\mathcal{O}_2}(\rho_0,\rho_1)$ then the interpolation $\rho_t:=((1-t)P_1+tP_2)_\#\Pi$
is absolutely continuous with respect to $\mu$ for any $t\in (0,1)$.

We are now able to pick up some elements in $\overline{\mathcal{O}_2}(\rho_0,\rho_1)$.

\begin{theorem}\label{NonEmptyO2}
	$\overline{\mathcal{O}_2}(\rho_0,\rho_1)$ is a non empty set.
\end{theorem}

\vskip 2mm
{\bf Proof.}
For all $n\in \mathbb{N}$ and $\ee>0$ we consider
$\Pi_{n,\ee} \in C(\rho_0^n,\rho_1^n)$ optimal for the cost $c_{\ee}$. It implies that
(\ref{ConvEnt1FiniCee0}) holds for $\Pi_{n,\ee}$. 
Now we pass to the Hilbert space and up to a subsequence, $(\Pi_{n,\ee})_n$ converges weakly to some coupling $\Pi_\ee \in C(\rho_0,\rho_1)$ which
solution of the problem $(P_\ee)$. Therefore (\ref{HConvEnt3Fini0}) holds for $\Pi_\ee$.
Again if $\Pi$ is a limit point of $(\Pi_\ee)_\ee$, then again (\ref{HConvEnt4Fini0}) holds for $\Pi$, namely
$\Pi$ satisfies the \emph{convexity property}.
We claim that any cluster point of $(\Pi_\ee)_\ee$ belongs to $\mathcal{O}_2(\rho_0,\rho_1)$. As a consequence,
the set $\overline{\mathcal{O}_2}(\rho_0,\rho_1)$ will be non empty.

Let $\Pi$ be a limit point of $(\Pi_\ee)_\ee$.\\
\noindent $\ast$ $\Pi \in \mathcal{O}_1(\rho_0,\rho_1)$. Indeed if $\Pi_0 \in \mathcal{O}_1(\rho_0,\rho_1)$, for 
 $\ee>0$: 
 \begin{eqnarray*}
 	\int |x-y| d\Pi_\ee &\leq & \int |x-y| d\Pi_\ee + \ee \int \alpha(x-y)d\Pi_\ee \\
 	&\leq & \int |x-y| d\Pi_0+\ee \int \alpha(x-y)d\Pi_0.
 \end{eqnarray*}
 Letting $\ee \rightarrow 0$,
 $$\int |x-y| d\Pi \leq \liminf_{\ee \rightarrow 0}\int |x-y| d\Pi_\ee \leq \int |x-y| d\Pi_0.$$
 
\noindent $\ast$ $\Pi \in \mathcal{O}_2(\rho_0,\rho_1)$. Indeed if $\Pi_0\in \mathcal{O}_2(\rho_0,\rho_1)$, for
$\ee>0$:
	\begin{eqnarray*}
	\int |x-y| d\Pi_\ee +\ee \int \alpha(x-y)d\Pi_\ee &\leq & \int |x-y| d\Pi_0
	+\ee \int \alpha(x-y)d\Pi_0\\
	&\leq & \int |x-y|d\Pi_\ee +\ee \int \alpha(x-y)d\Pi_0,
	\end{eqnarray*}
	the latter inequality is provided by the fact that $\Pi_0$ belongs in particular to $\mathcal{O}_1(\rho_0,\rho_1)$.
	Remove the same terms, dividing by $\ee$ and letting $\ee \rightarrow 0$,
	$$\int \alpha(x-y)d\Pi\leq \liminf_{\ee \rightarrow 0}\int \alpha(x-y)d\Pi_\ee \leq \int \alpha(x-y)d\Pi_0.$$
\fin

\vskip 3mm

Note also that for $\Pi_1$ and $\Pi_2$ are two coupling in $C(\rho_0,\rho_1)$ enjoying the \emph{convexity property},
every linear combination $(1-t)\Pi_1+t\Pi_2$ still enjoys the \emph{convexity property}. As a consequence 
$\overline{\mathcal{O}_2}(\rho_0,\rho_1)$ is a convex set.

\section{Recalls on optimal transportation theory}

We refer to \cite{VIL} or \cite{USERG} for proofs of results of this section.
We denote by $Supp(\Pi)$ the \emph{support of} $\Pi$, namely the smallest closed subset of $W\times W$
on which $\Pi$ is concentrated.

\begin{definition}
	Let $(X,\mu)$ and $(Y,\nu)$ be two Polish probability spaces and 
	$c:X \times Y \longrightarrow [0,\infty]$ be a measurable \emph{cost function}. We say that $\Pi \in C(\mu,\nu)$ is
	$c-$\emph{cyclically monotone} when for any $N\in \mathbb{N}$ and $(x_1,y_1),\dots,(x_N,y_N)\in Supp(\Pi)$,
	we have:
	$$\sum_{i=1}^Nc(x_i,y_i)\leq \sum_{i=1}^Nc(x_i,y_{i+1}),$$
	with $y_{N+1}:=y_1$.
\end{definition}

\begin{proposition}\label{CyclicallyMonotone}
	Let $(X,\mu)$ and $(Y,\nu)$ be two Polish probability spaces and 
	$c:X \times Y \longrightarrow [0,\infty]$ be a lower semi-continuous cost function. Then any optimal  
	coupling of
	$$\min_{\Pi\in C(\mu,\nu)}\int_{X\times Y}c(x,y)d\Pi(x,y)$$
	is $c-$cyclically monotone.
\end{proposition}

\begin{proposition}\label{Potentials}
	Let $\mu$ and $\nu$ be two probability measures on a Polish space $X$ and $c:X\times X\longrightarrow [0,\infty)$
	a cost function induced by the distance on $X$ i.e. $c(x,y)=d(x,y)$.
	If $\Pi$ is optimal for the Monge-Kantorovich problem between $\mu$ and $\nu$ with respect 
	to the cost $c$, then we can find
	a $\mu-$measurable $1-$Lipschitz map $u:X\longrightarrow X$ such that:
	\begin{eqnarray}
	\left\{ \begin{array}{llll}
				u(x)-u(y) &=& c(x,y) & \forall (x,y)\in Supp(\Pi)\\
				u(x)-u(y) &\leq & c(x,y) & \textrm{otherwise}\end{array} \right.
	\end{eqnarray}
\end{proposition}

Because in our case the cost we are interested in is a distance $|.|$ over $H$, we consider a map $u$ taken from this Proposition \ref{Potentials}.
It is worth to notice that the Problem (\ref{WMK2Infini}) is the same as the following 
$$\min_{\Pi \in C(\rho_0,\rho_1)} \int_{H \times H}\beta(x,y)d\Pi(x,y),$$
where the cost $\beta$ is defined by 
\begin{eqnarray}
	\beta(x,y):= \left\{ \begin{array}{ll}
				\alpha(x-y) & \textrm{if }u(x)-u(y)=\|x-y\|_\infty\\
				+\infty& \textrm{otherwise}\end{array} \right.
\end{eqnarray}

We complete this section with the following Lemma, which is proved in \cite{CHAMPS} and easily adaptable
in our setting.

Let $\rho_0$ and $\rho_1$ be two Borel probability measures on $W$.

\begin{lemma}\label{HSupOpt0}
	If $\Pi\in \mathcal{O}_2(\rho_0,\rho_1)$ then $\Pi$ is concentrated on some $\sigma-$compact set $\Gamma$ satisfying:
	\begin{equation}\label{HSupOpt}
		\forall (x,y),(x',y')\in \Gamma,~~~~x\in [x',y'] \Rightarrow (\nabla \alpha(y-x')-\nabla \alpha(y'-x),x-x')\geq 0.
	\end{equation}	
\end{lemma}

\vskip 2mm
{\bf Proof.}
	Since $\Pi$ is a solution of (\ref{MK1Hilbert}), there is a  Borel subset $\Gamma$ of $H\times H$ which
	is $|.|-$cyclically monotone. By inner regularity, up to remove a Borel set of zero measure,
	we can take $\Gamma$ $\sigma-$compact.
	According to Proposition \ref{Potentials}, we can find a potential $u:H	\longrightarrow H$ such that:
	$$\forall (x,y)\in \Gamma,~~~~u(x)-u(y)=|x-y|.$$
	Let $(x,y),(x',y')\in \Gamma$ such that $x\in [x',y']$. We have then:
	\begin{eqnarray*}
		u(x) &=& u(y) + |x-y|,\\
		u(x') &=& u(y') + |x'-y'|,
	\end{eqnarray*}
	and since $x\in [x',y']$, we also have:
	$$ |x'-y'| = |x-x'|+|x-y'|.$$
	Our potential $u$ is a $1-$Lipschitz map, so:
	$$u(x')=u(y')+|x-x'| +|x-y'| \geq u(x)+|x-x'| \geq u(x').$$
	This equality leads to:
	\begin{eqnarray*}
		u(x') &=& u(x) +|x-x'|= u(y) +|x-y| +|x-x'|\\
				&\geq& u(y) +|y-x'| \geq u(x').
	\end{eqnarray*}
	With the previous notation, it turns out that $\beta(x',y)=\alpha(x'-y)$ and $\beta(x,y')=\alpha(x-y')$.
	Moreover thanks to Proposition \ref{CyclicallyMonotone}, we also know that $\Pi$ is $\beta-$cyclically monotone
	hence by symmetry of $\alpha$:
	\begin{eqnarray*}
		\alpha(y-x)+\alpha(y'-x') \leq \alpha(y'-x)+\alpha(y-x').
	\end{eqnarray*}
	But by convexity of $\alpha$, we have:
	\begin{eqnarray*}
		\alpha(y-x)-\alpha(y-x') &\geq & \nabla \alpha(y-x').(x'-x),\\
		\alpha(y'-x)-\alpha(y'-x')&\leq & -\nabla \alpha(y'-x).(x-x').
	\end{eqnarray*}
	So combining these inequalities with the $\alpha-$monotonicity we get:
	$$(\nabla \alpha(y-x')-\nabla \alpha(y'-x),x-x')\geq 0.$$
\fin

\vskip 2mm

\begin{remark}
	As in \cite{CHAMPS2} the only reason to deal with $\sigma-$compact set $\Gamma$, is that the projection $P_1(\Gamma)$
	is also $\sigma-$compact, and in particular a Borel set.
\end{remark}

\section{Structure of the support of some element of $\mathcal{O}_2(\rho_0,\rho_1)$}

\vskip 2mm

Throughout this part, \emph{Differentiation theorem} \ref{LebPoints} is used many times.
We will present results in general framework.
We consider $\Pi\in C(\rho_0,\rho_1)$ and $\Gamma\subset W\times W$ a $\sigma-$compact set on
which $\Pi$ is concentrated. For all the sequel
we assume that $\rho_0=f\mu$ (the first measure has a density $f$ w.r.t. $\mu$).

Let us fix a sequence of positive number $(\delta_p)_p$ which tends to $0$ when $p$ goes to infinity.
\vskip 2mm

The following Lemma is a reinforcement of the one in \cite{CHAMP} (Lemma 3.3).

\begin{lemma}\label{SuppHilbLem1}
	Let $(y_n)_n$ be a dense sequence in $H$.
	Then we can find a Borel subset $D(\Gamma)$ on which $\Pi$ is still concentrated
	and such that for all $(x,y)\in D(\Gamma)$, $\forall r>0$, there exist $n,k\in\mathbb{N}$ satisfying
	$y \in B(y_n,\frac{1}{k+1})\subset B(y,r)$, $x\in Leb(f)\cap Leb(f_{n,k})$ and for all $p\in \mathbb{N}$:
	 $$\|f_{n,k}|_{B(x,\delta_p)}\|_{L^{\infty}}>0,$$
	where $f_{n,k}$ is the density of $(P_1)_\#\Pi_{|H \times \bar{B}(y_n,\frac{1}{k+1})}$.
\end{lemma}

\vskip 2mm
{\bf Proof.}
 Let $\delta=\delta_p>0$ be fixed. 
 We can find a recovering of $H$ with countably balls $(B(x_m^{(p)},\delta/2))_m$. 
 For any $(n,k)\in \mathbb{N}^2$
 we consider $f_{n,k}$ the density of the first marginal of the restriction of $\Pi$ to $H \times \bar{B}(y_n,\frac{1}{k+1})$ w.r.t. $\mu$.
 Fix $n,k\in \mathbb{N}$ and consider  
 $$D_{n,k}(\delta):=\left(\cup_{m\in \mathbb{N}}\{x\in B(x_m^{(p)},\delta/2),~\|f_{n,k}|_{B(x,\delta)}\|_{L^\infty}=0\}\right)
 \times \bar{B}(y_n,\frac{1}{k+1}).$$
 It turns out that
 $$\Pi(D_{n,k}(\delta))\leq \sum_{m\in \mathbb{N}}\int_{B(x_m^{(p)},\delta/2) \backslash \{\|f_{n,k}|_{B(x,\delta)}\|_{L^\infty}>0\}}f_{n,k}(x)d\mu(x)=0.$$
 Besides since $\rho_0<<\gamma$ it is straightforward to see that for 
 $C_{n,k}:=H\backslash (Leb(f)\cap Leb(f_{n,k}))\times H$,
 $$\Pi(C_{n,k})=\rho_0\left(H\backslash (Leb(f)\cap Leb(f_{n,k}))\right)=0.$$
 Therefore $\Pi$ is concentrated on the set $D_\delta(\Gamma):=\Gamma \backslash (\cup_{n,k} (D_{n,k}(\delta)\cup C_{n,k}))$.

 Since $(\delta_p)_p$ is a countably sequence, it follows $D(\Gamma):=\cap_p D_{\delta_p}(\Gamma)$
 has the desired properties. Indeed for any $\delta_p>0$ if $(x,y)\in D_{\delta_p}(\Gamma)$, 
 by density we can find $m,n,k\in \mathbb{N}$ such that
 $x\in B(x_m^{(p)},\delta_p/2), y\in B(y_n,1/(k+1))\subset B(y,r)$. The result ensues. 

\fin

\vskip 2mm

Notice that the previous result is quite general, because it is true for any coupling, not necessarly \emph{optimal}.

\vskip 2mm

\begin{definition}\ref{Gamma-1}
	Let $\Gamma$ be a $\sigma-$compact subset of $H\times H$. For $y\in \Omega$ and $r>0$ we define:
	$$\Gamma^{-1}(\bar{B}(y,r)):=P_1\left(\Gamma \cap (H \times \bar{B}(y,r))\right).$$
	An element $(x,y)$ of $\Gamma$ is called $\Gamma-$\emph{regular point} if $x$ is a Lebesgue point
	of $\Gamma^{-1}(\bar{B}(y,r))$ for any $r>0$.
\end{definition}	

It is worth to noting that from the definition (\ref{Gamma-1}), for all measurable subset $A$ of $W$:

$$\Pi(A\times \bar{B}(y,r))=\Pi\left(A\cap \Gamma^{-1}(\bar{B}(y,r))\times \bar{B}(y,r)\right).$$

\vskip 2mm

\begin{lemma}\label{SupHilbLem2}
	Under assumptions of Lemma \ref{SuppHilbLem1}, any element of $D(\Gamma)$ is a $\Gamma-$regular point, namely :
	$$(x,y)\in D(\Gamma) \Longrightarrow \lim_{\delta\rightarrow 0} \frac{\mu(\Gamma^{-1}(\bar{B}(y,r))\cap B(x,\delta))}{\gamma(B(x,\delta))}=1.$$
\end{lemma}

For the sequel, we introduce the following notation : if $\Gamma \subset H\times H$ then $T(\Gamma)=\left\{(1-t)x+ty,~
(x,y)\in \Gamma \right\}$. Since $\Gamma$ is $\sigma-$compact, $T(\Gamma)$ is $\sigma-$compact as well.

\begin{proposition}\label{HSupProp}
 Let $\rho_0, \rho_1 \in D(Ent_\mu)$, and $\Pi \in \overline{\mathcal{O}_2}(\rho_0,\rho_1)$ concentrated on
 a $\sigma-$compact set $\Gamma$.
 Then for all $(x,y_0),(x,y_1)$ belonging to the set $D(\Gamma)$ obtained in the Lemma \ref{SuppHilbLem1}, 
 with $y_0\neq y_1$ and $\forall r>0$ taken such that the closed balls centered at $y_0$ and $y_1$ with radius $r$ are disjoint, it holds:
 $$\gamma\left(T\left(\Gamma \cap (B(x,\delta_p)\times B(y_0,r))\right)\cap \Gamma^{-1}(\bar{B}(y_1,r))\cap  B(x,2\delta_p)\right) >0,$$
 $\forall p\in \mathbb{N}$ large enough.
\end{proposition}

\vskip 2mm
{\bf Proof.}
	Let $f$ be the density of $\rho_0$ w.r.t. $\mu$.
	Consider $\Pi \in \overline{\mathcal{O}_2}(\rho_0,\rho_1)$
	and let $(x,y_0), (x,y_1) \in D(\Gamma)$
	 such that $y_0\neq y_1$. We can assume that $x\neq y_0$.
	We fix $r>0$ for that $\bar{B}(y_0,r)\cap \bar{B}(y_1,r)=\emptyset$. 
	Thanks to the discussion above (Lemma \ref{SuppHilbLem1}), we introduce $n_0,n_1,k\in \mathbb{N}$ such that $B(y_{n_0},\frac{1}{k+1})\subset B(y_0,r)$, $B(y_{n_1},\frac{1}{k+1})\subset B(y_1,r)$. 
  Since $\delta_p$ decreases to $0$, we find $p\in \mathbb{N}$ large enough so that ${0<\delta=\delta_p<|x-y_0|+r}$,
  and
 \begin{equation}\label{gammaGamma}
 		\gamma\left(B(x,\delta)\cap \Gamma^{-1}(\bar{B}(y_0,r))\cap \Gamma^{-1}(\bar{B}(y_1,r))\right)>0.
 \end{equation}
 	This latter fact is possible thanks to the Proposition \ref{SupHilbLem2}.
  The corresponding densities given by Lemma \ref{SuppHilbertLem1} are denoted by $f_{n_0,k},~f_{n_1,k}$.

	Let us consider the Borel (up to a negligible set) set
\begin{eqnarray*}
	G_{x}:=\{z\in B(x,\delta),~f_{n_0,k}(z)>0,~f_{n_1,k}(z)>0\}.
\end{eqnarray*}
It turns out that $\mu(G_x)>0$. Indeed according to Lemma \ref{SuppHilbertLem1}:
	\begin{eqnarray*}
		\|f_{n_0,k}|_{B(x,\delta)}\|_{L^\infty}>0,\\
		\|f_{n_1,k}|_{B(x,\delta)}\|_{L^\infty}>0.
	\end{eqnarray*}
	Moreover we have 
	\begin{eqnarray*}
		\int_{\Gamma^{-1}(\bar{B}(y_0,r))\cap B(x,\delta)}f_{n_0,k}d\mu>0,\\
		\int_{\Gamma^{-1}(\bar{B}(y_1,r))\cap B(x,\delta)}f_{n_1,k}d\mu>0.
	\end{eqnarray*}
	The claim ensues thanks to (\ref{gammaGamma}).
		
\vskip 2mm
		
	Because $f_{n_1,k}$
	is the density of $(P_1)_\#\Pi_{|W\times \bar{B}(y_{n_1},\frac{1}{k+1})}$ we notice that:
	\begin{eqnarray*}	
		\Pi\left(G_x \times \bar{B}(y_{n_1},\frac{1}{k+1})\right) &=& \Pi\left(G_x\cap \Gamma^{-1}(\bar{B}(y_{n_1},\frac{1}{k+1})) \times  \bar{B}(y_{n_1},\frac{1}{k+1})\right) \\
		\mathrm{hence}~~\int_{G_x}f_{n_1,k}d\mu &=& \int_{G_x\cap \Gamma^{-1}(\bar{B}(y_{n_1},\frac{1}{k+1}))}f_{n_1,k}d\mu >0.
	\end{eqnarray*}
	It follows that 
	\begin{equation}\label{WInfiniPrem0}
		\gamma(G_x\cap \Gamma^{-1}(\bar{B}(y_1,r)))\geq \gamma\left(G_x\cap \Gamma^{-1}(\bar{B}(y_{n_1},\frac{1}{k_1+1}))\right)>0.
	\end{equation}

	Let $A(\delta):=B(x,2\delta)\cap \Gamma^{-1}(\bar{B}(y_1,r))\cap T\left(\Gamma\cap (B(x,\delta)\times B(y_0,r))\right)$.

	Consider the set $A_{x}:=G_{x}\times \bar{B}(y_{n_0},\frac{1}{k+1})$, and denote by
	$\Pi_{A_x}$ the restriction of $\Pi$ on $A_x$. We fix
	from now $t\in (0,\frac{\delta}{\|x-y_0\|_\infty+r})$ so that:
	if $z\in B(x,\delta)$ and $w\in B(y_0,r)$ then $(1-t)z+tw\in B(x,2\delta)$. Indeed
	\begin{eqnarray*}
		|(1-t)z+tw-x| &\leq & (1-t)|z-x|+t|w-x|\\
		&\leq & |z-x|+t(|w-y_0|+|y_0-x|)\\
		&< & \delta+\delta =2\delta.
	\end{eqnarray*}
	Therefore if we define $\rho_t^{A_x}:=((1-t)P_1+tP_2)_\#\Pi_{A_x}$, firstly we have:
	 $$	(P_1)_\#\Pi_{A_x}(G_x) \leq(P_1)_\#\Pi_{A_x}(B(x,\delta)) \leq \rho_t^{A_x}(B(x,2\delta))$$
	and thus:
	$$ (P_1)_\#\Pi_{A_x}(G_x\cap \Gamma^{-1}(\bar{B}(y_1,r))) \leq \rho_t^{A_x}(B(x,2\delta)\cap \Gamma^{-1}(\bar{B}(y_1,r))).$$
	Secondly thanks to (\ref{WInfiniPrem0}):
	\begin{eqnarray*}
		(P_1)_\#\Pi_{A_x}(G_x\cap \Gamma^{-1}(\bar{B}(y_1,r))) &=& \Pi\left(G_x\cap \Gamma^{-1}(\bar{B}(y_1,r))\times \bar{B}(y_{n_0},\frac{1}{k+1})\right)\\
			&=&\int_{G_x\cap \Gamma^{-1}(\bar{B}(y_1,r/2))}f_{n_0,k}d\gamma >0.
	\end{eqnarray*}
	And we deduce
	\begin{equation}\label{WInfiniPremInStSup}
		\rho_t^{A_x}(B(x,2\delta)\cap \Gamma^{-1}(\bar{B}(y_1,r)))>0.
	\end{equation}
	
	On the other hand, notice that $\rho_t^{A_x}$ 
	is concentrated on $T(\Gamma\cap (B(x,\delta)\times B(y_0,r))$ hence:
	\begin{eqnarray*}
	 ~&	\rho_t^{A_x}(B(x,2\delta)\cap \Gamma^{-1}(\bar{B}(y_1,r))) \\
	 =& \rho_t^{A_x}\left(B(x,2\delta)\cap T\left(\Gamma \cap (B(x,\delta)\times B(y_0,r))\right)\cap \Gamma^{-1}(\bar{B}(y_1,r))\right).
	\end{eqnarray*}
	
	Combining this latter fact with (\ref{WInfiniPremInStSup}), we get:
	\begin{eqnarray*}
		\rho_t^{A_x}(A(\delta))>0.
	\end{eqnarray*}
	And we know that $\rho_t^{A_x}$ inherits of the convexity property, so
	 is absolutely continuous w.r.t. $\gamma$. Hence it implies
	$\gamma(A(\delta))>0$.

\fin

\vskip 2mm

It is worth to notice that the Lebesgue differentiation theorem (Theorem \ref{LebPoints}) is
only used to get the positivity in (\ref{gammaGamma}). This is provided by the Proposition \ref{SupHilbLem2},
which needs this theorem. Without the theorem \ref{LebPoints}, 
the set considered in (\ref{gammaGamma}) is still non empty because
it containts $x$, but it can be of null measure.

\section{Proof of the main theorem and comments}

This section is devoted to prove Theorem \ref{MongeHilbert}.

\begin{theorem}
	Let $\rho_0, \rho_1 \in D(Ent_\mu)$. If (\ref{MongeHilbert}) is finite for some coupling, then any element of  	$\overline{\mathcal{O}_2}(\rho_0,\rho_1)$ is induced by
	a map $T$, and therefore $\overline{\mathcal{O}_2}(\rho_0,\rho_1)$
	is reduced to one element.
\end{theorem}

\vskip 2mm
{\bf Proof.}
	Let $\Pi \in \overline{\mathcal{O}_2}(\rho_0,\rho_1)$. In particular $\Pi \in \mathcal{O}_2(\rho_0,\rho_1)$
	and is concentrated on a $\sigma-$compact set $\Gamma$ satisfying (\ref{HSupOpt}).
	Furthermore Lemma \ref{SuppHilbLem1} provides us a $\sigma-$compact 
	set $D(\Gamma)$ on which $\Pi$ is still concentrated.
	We claim that $D(\Gamma)$ is contained in a graph of some Borel map.
	Let $(x_0,y_0)$ and $(x_0,y_1)$ in $D(\Gamma)$ and suppose that $y_0\neq y_1$. We can also assume $x_0\neq y_0$. 
 	By strict convexity of $\alpha$ we have:
 	$$((y_1-x_0)-(y_0-x_0),\nabla \alpha(y_1-x_0)-\nabla \alpha(y_0-x_0))>0.$$
 	Hence either $(y_1-x_0,\nabla \alpha(y_1-x_0)-\nabla \alpha(y_0-x_0))$ or
 	$(y_0-x_0,\nabla \alpha(y_0-x_0)-\nabla \alpha(y_1-x_0))$ is positive.
 	So without lost of generality we assume that:
	$$(\nabla \alpha(y_1-x_0)-\nabla \alpha(y_0-x_0),y_0-x_0) <0.$$
	By continuity of $\nabla \alpha$ we can find $r>0$ small enough so that:
	\begin{equation}\label{HMPIn1}
		\forall x,x'\in B(x_0,r),~\forall y'\in B(y_0,r),~\forall y\in B(y_1,r):
		~~(\nabla \alpha(y-x')-\nabla \alpha(y'-x),y'-x)<0.
	\end{equation}
	$r>0$ can be chosen so that the balls $\bar{B}(y_0,r)$ and $\bar{B}(y_1,r)$ are disjoint.
	
	Applying Proposition \ref{HSupProp} to $\left((x_0,y_0),(x_0,y_1)\right)$ we get:
	 $$\mu\left(T\left(\Gamma \cap (B(x_0,\delta_p)\times B(y_0,r/2))\right)\cap \Gamma^{-1}(\bar{B}(y_1,r/2))\cap B(x_0,2\delta_p)\right) >0,$$
	 $\forall p\in \mathbb{N}$ large enough.
	As a consequence we can find a $\delta=\delta_p \in (0,r/2)$	small enough in such a way that 
	there exist $(x',y') \in \Gamma \cap B(x_0,\delta)\times B(y_0,r/2)$
	and $x\in [x',y']\cap B(x_0,2\delta)$ and $y$ such that:
	$$(x,y)\in \Gamma \cap \left(([x',y']\cap B(x_0,2\delta))\times B(y_1,r)\right).$$
	Since $x\in [x',y']$, we have $x-x'=\frac{|x-x'|}{|y'-x|}(y'-x)$. So by (\ref{HSupOpt}), we have:
	$$(\nabla \alpha(y-x')-\nabla \alpha(y'-x),x-x')=
	\frac{|x-x'|}{|y'-x|}(\nabla \alpha(y-x')-\nabla \alpha(y'-x),y'-x)\geq 0,$$
	which contradicts (\ref{HMPIn1}). We obtain $y_1=y_0$.

	The unicity ensues from the convexity of $\overline{\mathcal{O}_2}(\rho_0,\rho_1)$, by the usual argument.

\fin

\vspace{5pt}

Let us make some comments.

\vspace{5pt}
	
	 We have proved that $\overline{\mathcal{O}_2}(\rho_0,\rho_1)$ is reduced to one element. 
	 However we do not know if $\mathcal{O}_2(\rho_0,\rho_1)$ has a	unique element.

\vspace{5pt}

 In \cite{CHAMP}, the authors do not require the absolute continuity of $\rho_t$ because 
the Lebesgue measure is doubling and invariant by translations. Thanks to that they can obtain good
bounds for $\rho_t$ (see Proposition 2.2 in \cite{CHAMP}).

\vspace{5pt}

 The fact that $\rho_1$ is absolutely continuous with respect to $\gamma$ is important for the section 2,
but we could hope it is possible to show the absolute continuity of interpolations $\rho_t$ ($t<1$)
without to pass by section 2. If it would be the case, the theorem \ref{MongeHilbert} would be true
for any probability measure $\rho_1$.

\vspace{5pt}

The strategy presented through the paper is general in the sense that the Hilbert norm $|.|$ 
could be replaced by any finite-valued norm $\|.\|$ on the Hilbert space $H$. 

\nocite{*}
\bibliographystyle{plain}
\bibliography{biblio}

\end{document}